\documentclass[11pt]{amsart}

\newtheorem{theorem}{Theorem}

\newtheorem{corollary}[theorem]{Corollary}

\newtheorem{lemma}[theorem]{Lemma}

\newtheorem{proposition}[theorem]{Proposition}
\newtheorem{remark}[theorem]{Remark}

\usepackage[active]{srcltx} 

\def\J#1#2#3{ \left\{ #1,#2,#3 \right\} }

\def\NN{{\mathbb{N}}}
\def\11{{1}}
\def\CC{{\mathbb{C}}}

\usepackage{enumerate}
\usepackage{amsmath}
\usepackage{amssymb}
\usepackage{hyperref}

\begin{document}

\title[Automatic continuity of biorthogonality preservers]{Automatic continuity of biorthogonality preservers
between compact C$^*$-algebras and von Neumann algebras}

\author[Burgos]{Mar{\' i}a Burgos}
\email{mariaburgos@ugr.es}
\address{Departamento de An{\'a}lisis Matem{\'a}tico, Facultad de
Ciencias, Universidad de Granada, 18071 Granada, Spain.}

\author[Garc{\' e}s]{Jorge J. Garc{\' e}s}
\email{jgarces@ugr.es}
\address{Departamento de An{\'a}lisis Matem{\'a}tico, Facultad de
Ciencias, Universidad de Granada, 18071 Granada, Spain.}

\author[Peralta]{Antonio M. Peralta}
\email{aperalta@ugr.es}
\address{Departamento de An{\'a}lisis Matem{\'a}tico, Facultad de
Ciencias, Universidad de Granada, 18071 Granada, Spain.}

\thanks{Authors partially supported by I+D MEC project no. MTM 2008-02186, and Junta de Andaluc{\'\i}a grants FQM 0199 and FQM 3737.}
\thanks{Published at Journal of Mathematical Analysis and Applications \url{https://doi.org/https://doi.org/10.1016/j.jmaa.2010.11.007}. This manuscript version is made available under the CC-BY-NC-ND 4.0 license \url{https://creativecommons.org/licenses/by-nc-nd/4.0/}}

\date{}

\begin{abstract} We prove that every biorthogonality preserving linear surjection between two
dual or compact C$^*$-algebras or between two von Neumann algebras is automatically continuous.
\end{abstract}

\maketitle
 \thispagestyle{empty}

\section{Introduction and preliminaries}

Two elements $a,b$ in a C$^*$-algebra $A$ are said to be
\emph{orthogonal} (denoted by $a \perp b$) if $a b^* = b^* a=0$. A
linear mapping $T: A\rightarrow B$ between two C$^*$-algebras is
called \emph{orthogonality preserving} if $T(x) \perp T(y)$
whenever $x \perp y$. The mapping $T$ is \emph{biorthogonality
preserving} whenever the equivalence $$ x \perp y \Leftrightarrow
T(x) \perp T(y)$$ holds for all $x,y$ in $A$.\smallskip

It can easily be seen that every biorthogonality preserving linear
surjection, $T:A\to B$ between two C$^*$-algebras is injective.
Indeed, for each $x\in A$, the condition $T(x)=0$ implies that
$T(x)\perp T(x)$, and hence $x\perp
 x$, which gives $x=0$.\smallskip

The study of orthogonality preserving operators between
C*-algebras begins with the work of W. Arendt \cite{Arendt} in the
setting of unital abelian C$^*$-algebras. His main result
establishes that every orthogonality preserving bounded linear
mapping $T : C(K) \to C(K)$ is of the form $$T(f) (t) = h(t)
f(\varphi (t)) \ (f\in C(K), \ t\in K),$$ where $h\in C(K)$ and
$\varphi: K \to K$ is a mapping which is continuous on $\{ t\in K
: h(t) \neq 0\}$. Several years later, K. Jarosz \cite{Jar}
extended the study to the setting of orthogonality preserving (not
necessarily bounded) linear mappings between abelian
C$^*$-algebras.\smallskip

A linear mapping $T: A\rightarrow B$ between two C$^*$-algebras is
said to be \emph{symmetric} if $T(x)^* = T(x^*)$, equivalently,
$T$ maps the self-adjoint part of $A$ into the self-adjoint part
of $B$.\smallskip

The study of orthogonality preservers between general
C$^*$-algebras was started in \cite{Wol}. M. Wolff proved in
\cite[Theorem 2.3]{Wol} that every orthogonality preserving
bounded linear and symmetric mapping between two C*-algebras is a
multiple of an appropriate Jordan $^*$-homomorphism.\smallskip

F.J. Fern\'andez-Polo, J. Mart\'inez and the authors of this note
studied and described orthogonality preserving bounded linear maps
between C$^*$-algebras, JB$^*$-algebras and JB$^*$-triples in
\cite{BurFerGarMarPe} and \cite{BurFerGarPe}. New techniques
developed in the setting of JB$^*$-algebras and JB$^*$-triples
were a fundamental tool to establish a complete description of all
orthogonality preserving bounded linear (non necessarily
symmetric) maps between two C*-algebras. We recall some background
results before stating the description obtained.\smallskip

Each C$^*$-algebra $A$ admits a Jordan triple product defined by
the expression $\J abc := \frac12 (a b^* c +c b^*a)$. Fixed points
of the triple product are called \emph{partial isometries} or
\emph{tripotents}. Every partial isometry $e$ in $A$ induces a
decomposition of $A$ in the form
$$A = A_2 (e) \oplus A_1(e) \oplus A_0 (e),$$ where $A_2 (e) := ee^* A e^*e$,
$A_1 (e) := (1-ee^*) A e^*e\oplus ee^* A (1-e^*e),$ and $A_0 (e)
:= (1-ee^*) A (1-e^*e)$. This decomposition is termed the
\emph{Peirce decomposition}. The subspace $A_2 (e) $ also admits a
structure of unital JB$^*$-algebra with product and involution
given by $x\circ_e y := \J xey$ and $x^{\sharp_e} := \J exe$,
respectively (compare \cite{HoS}). The element $e$ acts as the
unit element of $A_2 (e)$ (we refer to \cite{HoS} and \cite{Top}
for the basic results on JB- and JB$^*$-algebras).\smallskip

Recall that two elements $a$ and $b$ in a JB*-algebra
$J$ are said to \emph{operator commute} in $J$ if the
multiplication operators $M_a$ and $M_b$ commute, where $M_a$ is
defined by $M_a (x) := a\circ x$. That is, $a$ and $b$ operator
commute if and only if
 $(a\circ x) \circ b = a\circ (x\circ b)$ for
all $x$ in $J$. Self-adjoint elements $a$ and $b$ in $J$ generate
a JB$^*$-subalgebra that can be realised as a JC$^*$-subalgebra of
some $B(H)$ (c.f. \cite{Wri77}), and, in this realisation, $a$ and $b$
commute in the usual sense whenever they operator commute in $J$
\cite[Proposition 1]{Top}.  Similarly, two self-adjoint elements
$a$ and $b$ in $J$ operator commute if and only if $a^2 \circ b
=\J aba$ (i.e., $a^2 \circ b = 2 (a\circ b)\circ a - a^2 \circ
b$). If $b\in J$ we use $\{b\}^{'}$ to denote the set of elements
in $J$ that operator commute with $b$. (This corresponds to the
usual notation in von Neumann algebras).\smallskip

For each element $a$ in a von Neumann
algebra $W$ there exists a unique partial isometry $r(a)$ in $W$
such that $a=r(a) |a|$, and $r(a)^* r(a)$ is the support
projection of $|a|$, where $|a| =(a^* a)^{\frac12}$. It is also
known that $r(a) a^* r(a) = \J {r(a)}{a}{r(a)} = a$. We refer to
\cite[\S 1.12]{Sa} for a detailed proof of these results. The
element $r(a)$ will be called the \emph{range partial isometry}
of $a$.\smallskip

The characterisation of all orthogonality preserving bounded
linear maps between C*-algebras reads as follows:

\begin{theorem}\label{t BurFerGarMarPe}\cite[Theorem 17 and Corollary 18]{BurFerGarMarPe}  Let $T: A\to B$ be a bounded linear mapping between two C*-algebras. For $h= T^{**}
(1)$ and $r= r(h)$ the following assertions are equivalent.
\begin{enumerate}[{\rm $a)$}] \item $T$ is orthogonality
preserving.\item There exists a unique triple homomorphism $S: A
\to B^{**}$ satisfying $h^* S(z) = S(z^*)^* h,$ $h S(z^*)^* = S(z)
h^*,$ and $$T(z) = \frac12 \left( h r(h)^* S(z) + S(z) r(h)^*
h\right) $$ $$= h r(h)^* S(z)=S(z) r(h)^* h,$$ for all $z\in A$.

\item There exists a unique Jordan $^*$-homomorphism $S:
A \to B_2^{**} (r) $ satisfying that $S^{**} (1) = r$, $T(A)\subseteq \{h\}^{'}$ and
$T(z) = h\circ_r S(z) $ for all $z\in A$.

\item $T$ preserves zero triple products, that is, $\J {T(x)}{T(y)}{T(z)} = 0$ whenever $\J xyz =0$.$\hfill\Box$ \end{enumerate}
\end{theorem}

The problem of automatic continuity of those linear maps
preserving zero-products between C$^*$-algebras has inspired many
papers in the last twenty years.  A linear mapping between abelian
C$^*$-algebras is zero-products preserving if and only if it is
orthogonality preserving, however the equivalence doesn't hold for
general C$^*$-algebras (compare \cite[Comments before Corollary
18]{BurFerGarMarPe}). K. Jarosz proved the automatic continuity of
every linear bijection preserving zero-products between $C(K)$
spaces (see  \cite[Corollary]{Jar}). In 2003, M. A. Chebotar,
W.-F. Ke, P.-H. Lee, and N.-C. Wong showed that every
zero-products preserving linear bijection from a properly infinite
von Neumann algebra into a unital ring is automatically continuous
\cite[Theorem 4.2]{ChebKeLeeWo}. In the same year, J. Araujo and
K. Jarosz proved that a linear bijection which preserves
zero-products in both directions between algebras $L(X),$ of
continuous linear maps on a Banach space $X$, is automatically
continuous and a multiple of an algebra isomorphism
\cite{ArauJar}. The same authors conjectured that every linear
bijection between two C$^*$-algebras preserving zero-products in 
both directions is automatically continuous (see \cite[Conjecture
1]{ArauJar}).\smallskip

In this paper we study the problem of automatic continuity of
biorthogonality preserving linear surjections between
C$^*$-algebras. In Section 2 we prove that every biorthogonality
preserving linear surjection between two compact C$^*$-algebras is
continuous. In Sections 3 and 4, we establish, among many other
results, that every biorthogonality preserving linear surjection
between two von Neumann algebras is automatically continuous. It
follows as a consequence that a symmetric linear bijection between 
two von Neumann algebras preserving zero-products in both directions 
is automatically continuous. This provides a partial
answer to the conjecture posed by Araujo and Jarosz.

\subsection{Preliminary results}

A subspace $J$ of a C$^*$-algebra $A$ is said to be an \emph{inner
ideal} of A if $\{J, A, J\}\subseteq J$. Inner ideals in
C*-algebras were completely described by M. Edwards and F.
R\"{u}ttimann in \cite{EdRu91}.

Given a subset $M$ of $A,$ we write $M^\perp_{_A}$ for the
\emph{annihilator of $M$} (in $A$) defined by $$ M_{_A}^\perp:=\{
y \in A : y \perp x , \forall x \in M \}.$$ When no confusion can
arise, we shall write $M^{\perp}$ instead of $M^{\perp}_{_A}$. The
following properties can be easily checked.

\begin{lemma}
\label{l basic prop annihilator} Let $M$ be a subset of a
C$^*$-algebra $A$. The following statements
hold.\begin{enumerate}[{\rm $ a)$}] \item $M^{\perp}$ is a norm
closed inner ideal of $A$. When $A$ is a von Neumann algebra, then
$M^{\perp}$ is weak$^*$ closed. \item $M\cap M^{\perp}= \{0\},$
and $M\subseteq M^{\perp \perp}.$ \item $S^{\perp}\subseteq
M^{\perp}$ whenever $M\subseteq S\subseteq A.$ \item $M^{\perp}$
is closed for the product of $A$. \item $M^{\perp}$ is
$*$-invariant whenever $M$ is.$\hfill\Box$
\end{enumerate}
\end{lemma}

The next lemma describes the annihilator of a projection.

\begin{lemma}
\label{l bi annihilator closed projection} Let $p$ be a projection
in a (non necessarily unital) C$^*$-algebra $A$. The following assertions hold:
\begin{enumerate}[{\rm $ a)$}] \item $\{p\}_{_A}^{\perp} = (1-p) A (1-p)$,
where $1$ denotes the unit of $A^{**}$; \item $\{p\}_{_A}^{\perp\perp} = p A p.$
\end{enumerate}
\end{lemma}

\begin{proof} Statement {\rm $ a)$} follows straightforwardly.\smallskip

{\rm $ b)$} It is clear from {\rm $ a)$} that
$\{p\}_{_A}^{\perp\perp} \supseteq p A p.$ To to show the opposite
inclusion, let $a\in\{p\}_{_A}^{\perp\perp}$. Goldstine's theorem
(cf. Theorem V.4.2.5 in \cite{DunSchw}) guarantees that the closed
unit ball of $A$ is weak*-dense in the closed unit ball of
$A^{**}$. Thus, there exists a net $(x_{\lambda})$ in the closed
unit ball of $A$, converging in the weak*-topology of $A^{**}$ to
$1-p$. Noticing that $((1-p)x_{\lambda}(1-p))\subset
\{p\}_{_A}^{\perp}$, we deduce that \begin{equation} \label{eq
bi-anhi proj} (1-p)x_{\lambda}(1-p) a^* = a^*
(1-p)x_{\lambda}(1-p)=0,
\end{equation} for all $\lambda$.
Since the product of $A^{**}$ is separately weak*-continuous
(compare \cite[Theorem 1.7.8]{Sa}), taking weak*-limits in
(\ref{eq bi-anhi proj}), we have $(1-p) a^* = a^* (1-p) = 0$,
which shows that $p a = a p = a$.
\end{proof}

We shall also need some information about the norm closed inner
ideal generated by a single element. Let $a$ be an element in a
C$^*$-algebra $A$. Then $r(a) A^{**} r(a) \cap A= r(a) r(a)^*
A^{**} r(a)^* r(a) \cap A$ is the smallest norm closed inner ideal
in $A$ containing $a$ and will be denoted by $A(a)$. Further,  the
weak$^*$ closure of $A(a)$ coincides with $r(a) r(a)^* A^{**}
r(a)^* r(a)$ (c.f. \cite[Lemma 3.7 and Theorem 3.10 and its
proof]{EdRu91}). Since $\{a\}^{\perp\perp}$ is an inner ideal
containing $a$, we deduce that $A(a) \subseteq
\{a\}^{\perp\perp}$.\smallskip

It is well known that $\| \lambda a + \mu b\| = \max \{ \| \lambda
a \| , \|\mu b\|\}$, whenever $a\perp b$ and $\lambda,\mu\in
\mathbb{C}$. For every family $(A_i)_i$ of C$^*$-algebras, the
direct sum $\oplus^{\infty} A_i$ is another C$^*$-algebra with
respect to the pointwise product and involution. In this case, for
each $i\neq j$, $A_i $ and $A_j$ are mutually orthogonal
C$^*$-subalgebras of $\oplus_{i}^{\infty} A_i$.\smallskip

\begin{proposition}
\label{p direct sums} Let $A_1$, $A_2$ and $B$ be C$^*$-algebras
(respectively, von Neumann algebras). Let us suppose that $T :
A_1\oplus^{\infty} A_2 \to B$ is a biorthogonality preserving
linear surjection. Then $T(A_1)$ and $T(A_2)$ are norm closed
(respectively, weak$^*$ closed) inner ideals of $B$,
$B=T(A_1)\oplus^{\infty} T(A_2),$ and for $j= 1,2$ $T|_{A_j} : A_j
\to T(A_j)$ is a biorthogonality preserving linear surjection.
Further, if $T$ is symmetric then $T(A_1)$ and $T(A_2)$ are norm
closed (respectively, weak$^*$ closed) ideals of $B$.
\end{proposition}

\begin{proof}
Let us fix $j\in \{1,2\}$. Since $A_{j} = A_{j}^{\perp\perp}$ and
$T$ is a biorthogonality preserving linear surjection, we deduce
that $T(A_j) = T(A_j^{\perp\perp}) = T(A_j)^{\perp\perp}$. Lemma
\ref{l basic prop annihilator} guarantees that $T(A_j)$ is a norm
closed inner ideal of $B$ (respectively, a weak$^*$ closed
subalgebra of $B$ whenever $A_1$, $A_2$ and $B$ are von Neumann
algebras). The rest of the proof follows from Lemma \ref{l basic
prop annihilator} $e)$, and the fact that $B$ coincides with the
orthogonal sum of $T(A_1)$ and $T(A_2)$.
\end{proof}

\section{Biorthogonality preservers between dual C*-algebras}

A projection $p$ in a C$^*$-algebra $A$ is said to be
\emph{minimal} if $pAp = \mathbb{C} p.$ A partial isometry $e$ in
$A$ is said to be minimal if $e e^*$ (equivalently, $e^* e$) is a
minimal projection. The \emph{socle} of $A$, soc$(A)$, is defined
as the linear span of all minimal projections in $A$. The
\emph{ideal of compact elements} in $A$, $K(A)$, is defined as the norm
closure of soc$(A)$. A C$^*$-algebra is said to be \emph{dual} or
\emph{compact} if $A= K(A)$. We refer to \cite[\S 2]{Kap},
\cite{Alex} and \cite{Yli} for the basic references on dual
C$^*$-algebras.\smallskip

The following theorem proves that biorthogonality preserving
linear surjections between C$^*$-algebras send minimal projections to
scalar multiples of minimal partial isometries.

\begin{theorem}\label{thm minimal.to.minimal} Let $T: A\rightarrow B$ be a
biorthogonality preserving linear surjection between two
C$^*$-algebras and let $p$ be a minimal projection in $A$. Then
$\|T(p)\|^{-1} \ T (p) = e_p$ is a minimal partial isometry in
$B$. Further, $T$ satisfies that $T (p A p) = e_p e_p^* B e_p^*
e_p $ and $T((1-p) A (1-p)) = (1-e_p e_p^*) B (1-e_p^* e_p)$.
\end{theorem}

\begin{proof} Since $T$ is a biorthogonality preserving linear surjection, the equality
$$T(S_{_A}^\perp) = T(S)_{_B}^\perp$$ holds for every subset $S$ of
$A.$ For each minimal projection $p$ in $A$, $\{T(p)\}_{_B}^{
\perp \perp } = T(\{p\}_{_A}^{ \perp \perp })$ is an norm closed
inner ideal in $B$. Since $\{p\}_{_A}^{ \perp \perp } = p A p =
\CC p$, it follows that $\{T(p)\}_{_B}^{ \perp \perp }$ is a
one-dimensional subspace of $B$. Having in mind that
$\{T(p)\}_{_B}^{ \perp \perp }$ contains the inner ideal of $B$
generated by $T(p)$, we deduce that $B(T(p))$ must be
one-dimensional. In particular $(r(T(p)) r(T(p))^*) B^{**}
(r(T(p))^* r(T(p)) )= \overline{B(T(p))}^{w^*}$ has dimension one,
and hence $B(T(p))=(r(T(p)) r(T(p))^*) B^{**} (r(T(p))^* r(T(p)) )
= \CC r(T(p))$. This implies that $\|T(p)\|^{-1} \ T (p) = e_p$ is
a minimal partial isometry in $B$.\smallskip

The equality $T (p A p) = e_p e_p^* B e_p^* e_p $ has been proved.
Finally, $$T((1-p) A (1-p)) = T((p A p
)_{_{A}}^{\perp}) =\left( T(p A p)\right)_{_B}^{\perp} $$ $$= \left(e_p e_p^* B e_p^* e_p \right)_{_B}^{\perp} = (1-e_p e_p^*) B (1-e_p^* e_p).$$
\end{proof}

Let $a$ and $b$ be two elements in a C$^*$-algebra $A$. It is not hard to see that
$a\perp b$ if and only if $r(a)$ and $r(b)$ are two orthogonal partial isometries in $A^{**}$
(compare \cite[Lemma 1]{BurFerGarMarPe}).\smallskip

We shall make use of the following result which is a direct
consequence of Theorem \ref{t BurFerGarMarPe}. The proof is left for the reader.

\begin{corollary}\label{cor t BurFerGarMarPe} Let $T: A\to B$ be a bounded
 linear operator between two von Neumann algebras. For $h=
T (1)\in B$ and $r = r(h)$ the following assertions are
equivalent.
\begin{enumerate}[{\rm $a)$}] \item $T$ is a biorthogonality
preserving linear surjection. \item $h$ is invertible and there
exists a unique triple isomorphism $S: A \to B$ satisfying $h^*
S(z) = S(z^*)^* h,$ $h S(z^*)^* = S(z) h^*,$ and $$T(z) = \frac12
\left( h r(h)^* S(z) + S(z) r(h)^* h\right) $$ $$= h r(h)^*
S(z)=S(z) r(h)^* h,$$ for all $z\in A$.

\item $h$ is positive and invertible in $B_2 (r)$ and there exists
a unique Jordan $^*$-isomorphism $S: A \to B_2 (r) =B$ satisfying
that $S^{**} (1) = r$, $T(A)\subseteq \{h\}^{'}$ and $T(z) =
h\circ_r S(z) $ for all $z\in A$.
\end{enumerate} Further, in any of the previous statements, when $A$ is a factor, then $h$ is a multiple of the unit element in $B$.$\hfill\Box$
\end{corollary}

We deal now with dual C$^*$-algebras.

\begin{remark}
\label{r Cauchy series}{\rm Given a sequence $(\mu_n)\subset c_0$
and a bounded sequence $(x_n)$ in a Banach space $X$, the series
$\sum_{k} \mu_k x_k$ needs not be, in general, convergent in $X$.
However, when $(x_n)$ is a bounded sequence of mutually orthogonal
elements in a C$^*$-algebra, $A$, the equality
$$\left\| \sum_{k=1}^{n} \mu_k x_k - \sum_{k=1}^{m} \mu_k x_k \right\| =
\max \left\{ |\mu_{n+1}|,\ldots , |\mu_{m}|\right\}
\sup_{n+1 \leq k\leq m}\left\{\|x_k\|\right\},$$ holds for every $n<m$ in $\mathbb{N}$. It
follows that $\left(\sum_{k=1}^{n} \mu_k x_k\right)$ is a Cauchy series and
hence convergent in $A$. Alternatively, noticing that $\sum_{k}
x_k$ defines a \emph{w.u.C.} series in the terminology of
\cite{Die}, the final statement also follows from \cite[Theorem
V.6]{Die}.}
\end{remark}

\begin{lemma}\label{l c0 sums} Let $T: A\rightarrow B$ be a
biorthogonality preserving linear surjection between two
C$^*$-algebras and let $(p_n)_{n}$ be a sequence of
mutually orthogonal minimal projections in $A$.
Then the sequence $(\|T(p_n)\| )$ is bounded.
\end{lemma}

\begin{proof} By Theorem \ref{thm minimal.to.minimal}, for each natural $n$,
there exist a minimal partial isometry $e_n\in B$ and
$\lambda_n\in \mathbb{C}\setminus \{0\}$ such that $T(p_n) =
\lambda_n e_n$, and $\|T(e_n)\| =  \lambda_n$. Note that, by
hypothesis, $(e_n)$ is a sequence of mutually orthogonal minimal
partial isometries in $B$.\smallskip

Let $(\mu_n)$ be any sequence in $c_0$. Since the $p_n$'s are
mutually orthogonal, the series $\sum_{k\geq 1}  \mu_k p_k$
converges to an element in $A$ (compare Remark \ref{r Cauchy
series}). For each natural $n$, $\sum_{k\geq 1}^{\infty}  \mu_k
p_k $ decomposes as the orthogonal sum of $\mu_n p_n$ and
$\sum_{k\neq n}^{\infty}  \mu_k p_k $, therefore
$$T\left(\sum_{k\geq 1}^{\infty}  \mu_k p_k \right) = \mu_n
\lambda_n e_n  + T\left(\sum_{k\neq n}^{\infty}  \mu_k
p_k\right),$$ with $\mu_n \lambda_n e_n  \perp T\left(\sum_{k\neq
n}^{\infty}  \mu_k p_k\right)$, which in particular implies
$$\left\|T\left(\sum_{k\geq 1}^{\infty}  \mu_k p_k \right)
\right\|= \max\left\{|\mu_n| | \lambda_n|, \left\|
T\left(\sum_{k\neq n}^{\infty}  \mu_k
p_k\right)\right\|\right\}\geq |\mu_n| | \lambda_n|. $$ This
establishes that for each $(\mu_n)$ in $c_0$, $(\mu_n  \lambda_n)$
is a bounded sequence, which proves the statement.
\end{proof}

\begin{lemma}\label{l lim series to zero}
Let $T: A\rightarrow B$ be a biorthogonality preserving linear
surjection between two C$^*$-algebras, $(\mu_n)$ a sequence in
$c_0$ and let $(p_n)_{n}$ be a sequence of mutually orthogonal
minimal projections in $A$. Then the sequence $\left( T \left(
\sum_{k\geq n}^{\infty}  \mu_k p_k \right)\right)_{n}$ is well
defined and converges in norm to zero.
\end{lemma}

\begin{proof} By Theorem \ref{thm minimal.to.minimal} and
Lemma \ref{l c0 sums} it follows that $(T(p_n))$ is a bounded
sequence of mutually orthogonal elements in $B$. Let $M= \sup \{
\|T(p_n)\|: n\in\NN\}$. For each natural $n$, Remark \ref{r Cauchy
series} assures that the series $\sum_{k\geq n}^{\infty}  \mu_k
p_k$ converges.\smallskip

Let us define $y_n := T \left( \sum_{k\geq n}^{\infty}  \mu_k p_k
\right)$. We claim that $(y_n)$ is a Cauchy sequence in $B$.
Indeed, given $n< m$ in $\mathbb{N}$, we have \begin{equation}
\label{eq lemma lim series to zero} \| y_n -y_m\| = \left\| T
\left( \sum_{k\geq n}^{m-1}  \mu_k p_k \right) \right\| = \left\|
\sum_{k\geq n}^{m-1}  \mu_k T(p_k) \right\|
\end{equation} $$ \stackrel{(*)}{\leq} M\ \max \{ |\mu_{n}|,\ldots , |\mu_{m-1}|\},$$
where at $(*)$ we apply the fact that $(T(p_n))$ is a sequence of
mutually orthogonal elements. Consequently, $(y_n)$ converges in
norm to some element $y_0$ in $B$. Let $z_0$ denote $T^{-1}
(y_0)$.\smallskip

Let us fix a natural $m$. By hypothesis, for each $n>m$, $p_m$ is
orthogonal to $\sum_{k\geq n}^{\infty}  \mu_k p_k$. This implies
that $T(p_m) \perp y_n$, for every $n>m$, which, in particular,
gives $T(p_m)^* y_n = y_n T(p_m)^* =0$, for every $n>m$. Taking
limits when $n$ tends to $\infty$ we have $T(p_m)^* y_0 = y_0
T(p_m)^* =0$. This shows that $y_0=T(z_0) \perp T(p_m)$, and hence
$p_m\perp z_0$. Since $m$ was arbitrarily chosen we deduce that,
for each natural $n$, $z_0$ is orthogonal to $\sum_{k\geq
n}^{\infty}  \mu_k p_k$. Therefore, $(y_n)\subset \{
y_0\}_{_B}^{\perp},$ and hence $y_0 $ belongs to the norm closure
of $\{ y_0\}_{_B}^{\perp},$ which implies that $y_0 =0$.
\end{proof}

\begin{proposition}
\label{p weakly compact} Let $T: A\rightarrow B$ be a
biorthogonality preserving linear surjection between
C$^*$-algebras. Then $T|_{_{K(A)}}$ is continuous if and only if
the set $\left\{ \|T(p)\| : p \hbox{ minimal projection in }
A\right\} $ is bounded.
\end{proposition}

\begin{proof}
The necessity being obvious. Suppose that
$$M= \sup \left\{ \|T(p)\| : p \hbox{ minimal projection in } A\right\} <\infty.$$
Each nonzero self-adjoint element $x$ in $K(A)$ can be written as
a norm convergent (possibly finite) sum $x = \sum_{n} \lambda_n
p_n$, where $p_n$ are mutually orthogonal minimal projections in
$A$, and $\|x\| = \sup \{|\lambda_n| : n\}$ (compare \cite{Alex}).
If the series $x= \sum_{n} \lambda_n p_n$ is finite then $$\|T(x)
\| = \left\| \sum_{n=1}^{m} \lambda_n T(p_n) \right\|
\stackrel{(*)}{=} \max \left\{  \left\| \lambda_n T(p_n) \right\|
: n= 1,\ldots,m\right\} \leq M \|x\|,$$ where at $(*)$ we apply
the fact that $(T(p_n))$ is a finite set of mutually orthogonal
elements in $B$. When the series $x=\sum_{n} \lambda_n u_n$ is
infinite we may assume that $(\lambda_n)\in c_0$.\smallskip

It follows from Lemma \ref{l lim series to zero} that the sequence
$\left( T \left( \sum_{k\geq n}^{\infty}  \lambda_k p_k
\right)\right)_{n}$ is well defined and converges in norm to zero.
We can find a natural $m$ such that $\left\|T \left( \sum_{k\geq
m}^{\infty}  \lambda_k p_k \right)\right\| < M \|x\|$. Since the
elements $\lambda_1 p_1,\ldots,\lambda_{m-1} p_{m-1},$ $
\sum_{k\geq m}^{\infty}  \lambda_k p_k$ are mutually orthogonal,
we have $$\|T(x) \| = \max \big\{ \|T(\lambda_1
p_1)\|,\ldots,\|T(\lambda_{m-1} p_{m-1})\|, \left\|T(\sum_{k\geq
m}^{\infty}  \lambda_k p_k)\right\| \big\} \leq M \|x\|.$$ We have
established that $\|T(x)\| \leq M \|x\|$, for all $x\in
K(A)_{sa}$, and by linearity $\|T(x)\| \leq 2 M \|x\|$, for all
$x\in K(A).$
\end{proof}

\begin{theorem}
\label{t authom cont dual C-alg} Let $T: A\rightarrow B$ be a
biorthogonality preserving linear surjection between two
C$^*$-algebras. Then $T|_{_{K(A)}} : K(A) \to K(B)$ is continuous.
\end{theorem}

\begin{proof}
Theorem \ref{thm minimal.to.minimal} implies $T(\hbox{soc} (A))=
\hbox{soc} (B)$ (compare \cite[Theorem 5.1]{Yli}).  By Proposition
\ref{p weakly compact} it is enough to show the boundedness of the
set
$$\mathcal{P} = \left\{ \|T(p)\| : p \hbox{ minimal projection in
} A\right\} .$$

Suppose, on the contrary, that $\mathcal{P}$ is unbounded. We
shall show by induction that there exists a sequence $(p_n)$ of
mutually orthogonal minimal projections in $A$ such that
$\|T(p_n)\| > n$.\smallskip

The case $n=1$ is clear. The induction hypothesis guarantees the
existence of mutually orthogonal minimal projections
$p_1,\ldots,p_n$ in $A$ with $\|T(p_k)\| >k$ for all $k\in
\{1,\ldots,n\}$.\smallskip

By assumption, there exists a minimal projection $q\in A$
satisfying $\|T(q) \| > \max \{
\|T(p_1)\|,\ldots,\|T(p_n)\|,n+1\}$. We claim that $q$ must be
orthogonal to each $p_j$. Suppose, on the contrary, that for some
$j$, $p_j$ and $q$ are not orthogonal. Let $C$ denote the
C*-subalgebra of $A$ generated by $q$ and $p_j$. We conclude from
Theorem 1.3 in \cite{RaSin} (see also \cite[\S 3]{Ped68}) that
there exist $0<t<1$ and a $^*$-isomorphism $\Phi : C \to M_2
(\CC)$ such that $\Phi (p_j) = \left( \begin{array}{cc}
                                                     1 & 0 \\
                                                     0 & 0 \\
                                                   \end{array}
                                                 \right)$ and  $\Phi (q) = \left(
                                                                             \begin{array}{cc}
                                                                               t & \sqrt{t(1-t)} \\
                                                                               \sqrt{t(1-t)} & 1-t \\
                                                                             \end{array}
                                                                           \right).$
Since $T|_{C} :  C \cong M_2(\CC) \to T(C)$ is a continuous
biorthogonality preserving linear bijection, Theorem \ref{t
BurFerGarMarPe} (see also Corollary \ref{cor t BurFerGarMarPe})
guarantees the existence of an scalar $\lambda\in
\CC\backslash\{0\}$ and a triple isomorphism $\Psi : C \to T(C)$
such that $T(x) = \lambda \Psi (x)$ for all $x\in C$. In this
case, $\|T(p_j) \| = |\lambda| \|\Psi (p_j)\|$ implies that
$$\|T(p_j)\| < \|T(q)\| = |\lambda| \  \|\Psi (q)\| = |\lambda| \ \|\Psi (p_j)\| = \|T(p_j)\|,$$
which is a contradiction. Therefore $q\perp p_j$, for every $j=
1,\ldots,n$.\smallskip

It follows by induction that there exists a sequence  $(p_n)$ of
mutually orthogonal minimal projections in $A$ such that
$\|T(p_n)\| > n$. The series $\sum_{n=1}^{\infty}
\frac{1}{\sqrt{n}} p_n$ defines an element $a$ in $A$ (compare
Remark \ref{r Cauchy series}). For each natural $m$, $a$
decomposes as the orthogonal sum of $\frac{1}{\sqrt{m}} p_m$ and
$\sum_{n\neq m}^{\infty} \frac{1}{\sqrt{n}} p_n$, therefore
$$T(a)= \frac{1}{\sqrt{m}} T(p_m) + T\left(\sum_{n\neq m}^{\infty}
\frac{1}{\sqrt{n}} p_n\right),$$ with $\frac{1}{\sqrt{m}} T(p_m)
\perp T\left(\sum_{n\neq m}^{\infty} \frac{1}{\sqrt{n}}
p_n\right).$ This argument implies that $$\|T(a)\|=\max \left\{
\frac{1}{\sqrt{m}} \left\|T(p_m) \right\| ,\
\left\|T\left(\sum_{n\neq m}^{\infty} \frac{1}{\sqrt{n}}
p_n\right)\right\| \right\}> \sqrt{m}.$$ Since $m$ was arbitrarily
chosen, we have arrived at our desired contradiction.
\end{proof}

The following result is an immediate consequence of the above
theorem.

\begin{corollary}
\label{cor bop dual C*alg} Let $T: A\rightarrow B$ be a
biorthogonality preserving linear surjection between two
dual C$^*$-algebras. Then $T$ is continuous.
\end{corollary}

Given a complex Hilbert space $H$, it is well known that
soc$(L(H))$ coincides with the space of all finite rank operators
on $H$. The ideal $K(L(H))$ agrees with the ideal $K(H)$ of all
compact operators on $H$.

\begin{corollary}
\label{cor bop KH} Let $T: K(H)\rightarrow K(H)$ be a
biorthogonality preserving linear surjection, where $H$ is a complex Hilbert space.
Then $T$ is continuous.
\end{corollary}

\section{C$^*$-algebras linearly spanned by their projections}

In a large number of C$^*$-algebras every element can be expressed as a finite linear combination of projections:
the Bunce-Deddens algebras; the irrotational rotation algebras; simple, unital AF C*-algebras with finitely many
extremal states; UHF C*-algebras; unital, simple C$^*$-algebras of real rank zero with no tracial states; properly infinite C$^*$-
and von Neumann algebras; von Neumann algebras of type II$_1$ $\ldots$
(see \cite{Mar98}, \cite{MarMur02}, \cite{Mar06}, \cite{PearTopp}, \cite{LinMat} and the references therein).\smallskip

\begin{theorem}\label{expanded-by-proj}Let $T:A\to B$ be an orthogonality
preserving linear map between C$^*$-algebras, where $A$ is unital.
Suppose that every element of $A$ is a finite linear combination
of projections, then $T$ is continuous.
\end{theorem}
\begin{proof} Let $p$ be a projection in $A$.
As $p\perp (\11-p)$ then $T(p)\perp T(\11)-T(p)$,  that is
$T(p)T(\11)^*=T(p)T(p)^*$ and $T(\11)^*T(p)=T(p)^*T(p)$. In
particular, $T(p)T(\11)^*=T(\11)T(p)^*$ and
$T(\11)^*T(p)=T(p)^*T(\11)$. Since every element in $A$ coincides
with a finite linear combination of projections, it follows that
\begin{equation}\label{1}
T(x)T(\11)^*=T(\11)T(x^*)^*,
\end{equation}
for all $x\in A$.

Let now $p,q$ be two projections in $A$. The relation $qp \perp
(\11-q)(\11-p)$ implies that $T(qp)\perp T(\11-q-p+qp)$. Therefore
\begin{equation}\label{4}
T(qp)T(\11)^*-T(qp)T(q)^*-T(qp)T(p)^*+T(qp)T(qp)^*=0.
\end{equation}
Similarly, since $q(\11-p)\perp (\11-q)p$, we have $T(q-qp)\perp T(p-qp)$, and hence
\begin{equation}\label{6}
T(q)T(p)^*-T(q)T(qp)^*-T(qp)T(p)^*+T(qp)T(qp)^*=0.
\end{equation}

From (\ref{4}) and (\ref{6}), we get
$$T(qp)T(\11)^*-T(qp)T(q)^*=T(q)T(p)^*-T(q)T(qp)^*.$$

Being $A$ linearly spanned by its projections, the last equation
yields to
\begin{equation}\label{10}
T(qx)T(\11)^*-T(qx)T(q)^*=T(q)T(x^*)^*-T(q)T(qx^*)^*.
\end{equation} 
for all $x\in A$, and $q=q^*=q^2 \in A$.

By replacing, in (\ref{10}), $q$ with $\11-q$, we get
$$T(q-qx)T(\11)^*-T(x-qx)T(\11-q)^*=T(\11-q)T(x^*)^*-T(\11-q)T(x^*-qx^*)^*.$$
Having in mind (\ref{1}) we obtain
\begin{equation}\label{12}
T(x)T(q)^*-T(qx)T(q)^*=T(\11)T(qx^*)^*-T(q)T(qx^*)^*.
\end{equation}

From equations (\ref{10}) and (\ref{12}),
we deduce that
$$T(qx)T(\11)^*-T(x)T(q)^*=T(q)T(x^*)^*-T(\11)T(qx^*)^*.$$
for every $x$ in $A$ and every projection $q$ in $A$. Again,
the last equation and the hypothesis on $A$ prove:
$$T(yx)T(\11)^*-T(x)T(y^*)^*=T(y)T(x^*)^*-T(\11)T(y^*x^*)^*,$$  for all $x,y\in A$.
Since, by (\ref{1}), $T(\11)T(y^*x^*)^*=T(xy)T(\11)^*$, we can write the above equation as:
\begin{equation}\label{18}
T(yx+xy)T(\11)^*=T(y)T(x^*)^*+T(x)T(y^*)^*,
\end{equation} for all $x,y \in A$.\smallskip

Let $h=T(\11)$. It follows from (\ref{18}) that
$$T(x^2)h^*=T(x)T(x^*)^* \qquad(x\in A).$$ We claim that the
linear mapping $S:A\to B$, $S(x):=T(x)h^*$, is positive, and hence
continuous. Indeed, given $a\in A^+$, there exists $x\in A_{sa}$
such that $a=x^2$. Then $S(a)=T(x^2)h^*=T(x)T(x)^*\geq
0$.\smallskip

Finally, for any $x\in A_{sa}$, $$\|T(x)\|^2=\|T(x)T(x)^*\|=\|S(x^2)\|\leq \|S\| \|x\|^2,$$
which implies that $T$ is bounded on self-adjoint elements, and thus $T$ is continuous.
\end{proof}

We have actually proved the following:

\begin{proposition}\label{p expanded-by-proj}Let $T:A\to B$ be a
linear map between C$^*$-algebras, where $A$ is unital and every
element in $A$ is a finite linear combination of projections.
Suppose that $T$ satisfies one of the following statements:
\begin{enumerate}[$a)$] \item $a b^* =0 \Rightarrow T(a) T(b)^*
=0$. \item $ b^* a =0 \Rightarrow T(b)^*  T(a)=0$.
\end{enumerate} Then $T$ is continuous.$\hfill\Box$
\end{proposition}

Recall that a unital C$^*$-algebra is \emph{properly infinite} if
it contains two orthogonal projections equivalent to the identity
(i.e. it contains two isometries with mutually orthogonal range
projections). Zero product preserving linear mappings from a
properly infinite von Neumann algebra to a unital ring were
studied and described in \cite[Theorem 4.2]{ChebKeLeeWo}. In this
paper we consider a wider class of C$^*$-algebras. Let $A$ be a
properly infinite C$^*$-algebra or a von Neumann algebra of type
II$_1$. It follows by \cite[Corollary 2.2]{LinMat} (see also
\cite{PearTopp}) and \cite[Theorem 2.2.(a)]{GoldPas} that every
element in $A$ can be expressed as a finite linear combination of
projections. Our next result follows immediately from Theorem
\ref{expanded-by-proj}.

\begin{corollary}\label{prop-infinite} Let $A$ be a properly infinite
unital C$^*$-algebra or a von Neumann algebra of type II$_1$.
Every orthogonality preserving linear map from $A$ to another
C$^*$-algebra is automatically continuous.$\hfill\Box$
\end{corollary}

The following corollary is, in some sense, a generalisation of
\cite[Theorem 2]{Mol} and \cite[Theorem 1]{BaiHou}.

\begin{corollary}\label{c expanded-by-proj}Let $A$ be a properly infinite
unital C$^*$-algebra or a von Neumann algebra of type II$_1$ and
let $T:A\to B$ be a linear map from $A$ to another C$^*$-algebra.
Suppose that $T$ satisfies one of the following statements:
\begin{enumerate}[$a)$] \item $a b^* =0 \Rightarrow T(a) T(b)^*
=0$. \item $ b^* a =0 \Rightarrow T(b)^*  T(a)=0$.
\end{enumerate} Then $T$ is continuous.$\hfill\Box$
\end{corollary}

\section{Biorthogonality preservers between von Neumann algebras}

Let us recall some fundamental results derived from the Murray-von
Neumann dimension theory. Two projections $p$ and $q$ in a von
Neumann algebra $A$ are (Murray-von Neumann) \emph{equivalent}
(written $p\sim q$) if there exists a partial isometry $u\in A$
with $u^*u = p$ and $u u^* = q$. We write $q \precsim p$ when
$q\leq p$ and $p\sim q$. A projection $p$ in $A$ is said to be
\emph{finite} if $q \precsim p$ implies $p = q$. Otherwise, it is
called \emph{infinite}. A von Neumann algebra is said to be finite
or infinite according to the property of its identity projection.
A projection $p$ in $A$ is called \emph{abelian} if $p A p$ is a
commutative von Neumann algebra (compare \cite[\S
V.1]{Tak}).\smallskip

A von Neumann algebra $A$ is said to be of \emph{type} I if every
nonzero central projection in $A$ majorizes a nonzero abelian
projection. If there is no nonzero finite projection in $A$, that
is, if $A$ is purely infinite, then it is of \emph{type} III. If
$A$ has no nonzero abelian projection and if every nonzero central
projection in $A$ majorizes a nonzero finite projection of $A$,
then it is of \emph{type} II. If $A$ is finite and of type II
(respectively, type I), then it is said to be of type II$_1$
(respectively, type I$_{fin}$). If $A$ is of type II and has no
nonzero central finite projection, then $A$ is said to be of type
II$_{\infty}$. Every von Neumann algebra is uniquely decomposable
into the direct (orthogonal) sum of weak$^*$ closed ideals of type
I, type II$_1$, type II$_{\infty}$, and type III (this
decomposition is usually called, the \emph{Murray-von Neumann}
decomposition).\smallskip

\begin{proposition}
\label{p finite type I} Let $T: A \to B$ be a surjective linear
mapping from a unital C$^*$-algebra onto a finite von Neumann
algebra. Suppose that for each invertible element $x$ in $A$ we
have $\{T(x)\}^{\perp} \subseteq T(\{x\}^{\perp})$. Then $T$ is
continuous. Further, if $A$ is a von Neumann algebra, then $T (x)=
T(1) S(x)$ ($x\in A$), where $S: A \to B$ is a Jordan
homomorphism. In particular, every biorthogonality preserving
linear surjection between two von Neumann algebras one of which is
finite is continuous.
\end{proposition}

\begin{proof}
Let $T: A\to B$ be a surjective linear mapping satisfying the
hypothesis. We claim that $T$ preserves invertibility. Let $z$ be
an invertible element in $A$ and let $r= r(T(z))$ denote the range
partial isometry of $T(z)$ in $B$. In this case $$(1-r r^*) B
(1-r^* r) = \{T(z)\}^{\perp} \subseteq T(\{z\}^{\perp})= \{0\}. $$
Since $B$ is a finite von Neumann algebra, $1-r r^*$ and $1-r^* r$
are equivalent projections in $B$ (compare \cite[Exercise
6.9.6]{KadRing}). Thus, there exists a partial isometry $w$ in $B$
such that $w w^* = 1- r r^*$ and $w^* w= 1-r^* r$, and hence $ww^*
B w^* w=\{0\}$. It follows that $1- r r^*= ww^* = 0= w^* w= 1-r^*
r$.\smallskip

Since $1= r^* r$ (respectively, $1= r r^*$) is the support
projection of $T(z)^* T(z)$ (respectively, $T(z) T(z)^*$), we
deduce that $T(z) T(z)^*$ and $T(z)^* T(z)$ are invertible
elements in $B$, and therefore $T(z)$ is invertible.\smallskip

Having in mind that $T$ sends invertible elements to invertible
elements, we deduce from \cite[Corollary 2.4]{CuiHou} (see also
\cite[Theorem 5.5.2]{Aup}) that $T$ is continuous. If $A$ is a von
Neumann algebra, then it is well known that $T$ is a Jordan
homomorphism multiplied by $T(1)$ (c.f. \cite[Theorem 1.3]{Aup00}
or \cite[Corollary 2.4]{CuiHou}).
\end{proof}

It is obvious that every $^*$-isomorphism between two von Neumann
algebras preserves the summands appearing in the Murray-von
Neumann decomposition. However, it is not so clear that every
Jordan $^*$-isomorphism between two von Neumann algebras also
preserves the Murray-von Neumann decomposition. The justification
follows from an important result due to R. Kadison \cite{Kad}. If
$T: A \to B$ is a Jordan $^*$-isomorphism between von Neumann
algebras, then there exist weak$^*$ closed ideals $A_1$ and $A_2$
in $A$ and $B_1$ and $B_2$ in $B$ satisfying that $A = A_1
\oplus^{\infty} A_2$, $B= B_1 \oplus^{\infty} B_2$, $T|_{A_1} :
A_1 \to B_1$ is an $^*$-isomorphism, and $T|_{A_2} : A_2 \to B_2$
is an $^*$-anti-isomorphism (see \cite[Theorem 10]{Kad}). It
follows that every Jordan $^*$-isomorphism preserves the
Murray-von Neumann decomposition.

\begin{theorem}\label{t vonNeumann} Every biorthogonality preserving linear surjection between
von Neumann algebras is automatically continuous.
\end{theorem}

\begin{proof} Let $T : A\to B$ be a biorthogonality preserving linear surjection between
 von Neumann algebras.\smallskip

It is well known that every von Neumann algebra is uniquely decomposed into a
direct sum of five algebras, respectively, of types I$_{fin}$,
I$_{\infty}$, II$_1$, II$_{\infty}$ and III, where I$_{fin}$ is a
finite type I von Neumann algebra, II$_{1}$ is a finite type II von
Neumann algebra and the direct sum of those summands of types
I$_{\infty}$, II$_{\infty}$ and III is a properly infinite von
Neumann algebra (compare \cite[Theorem V.1.19]{Tak}). Therefore,
$A$ and $B$ decompose in the form $A = A_{I_{fin}}\oplus^{\infty}
A_{II_1}\oplus^{\infty} A_{p\infty}$, $B = B_{I_{fin}}\oplus^{\infty} B_{II_1}\oplus^{\infty}
B_{p\infty},$ where $ A_{I_{fin}}$ and $B_{I_{fin}}$ are finite
type I von Neumann algebras, $A_{II_1}$ and $B_{II_1}$ are type
II$_{1}$ von Neumann algebras, and $A_{p\infty}$ and
$B_{p\infty}$ are properly infinite von Neumann
algebras.\smallskip

Corollary \ref{prop-infinite} guarantees that $T|_{A_{p\infty}} :
A_{p\infty} \to B$ and $T|_{A_{II_1}} : A_{II_1} \to B$ are
continuous linear mappings. In order to simplify notation we
denote $A_1 = A_{I_{fin}}$, $A_2 = A_{II_1}\oplus^{\infty}
A_{p\infty}$, $B_1 = B_{I_{fin}}$, and $B_2 =
B_{II_1}\oplus^{\infty} B_{p\infty}$. According to this notation,
$T|_{A_2} : A_2 \to B$ is continuous. Theorem \ref{t
BurFerGarMarPe} assures the existence of a Jordan
$^*$-homomorphism $S_2:A_2 \to B_2 (r_2) $ satisfying that $S_2
(1_2) = r_2$, $T(A_2)\subseteq \{h_2\}^{'}$ and $T(z) =
h_2\circ_{r_2} S_2(z) $ for all $z\in A_2$, where $1_2$ is the
unit of $A_2$ and $r_2$ is the range partial isometry of
$T(1_2)=h_2$. We notice that, for each $z\in T(A_2)$, $r_2 r_2^* z
+ z r_2^* r_2 = 2 z,$ and we therefore have $r_2 r_2^* z = z r_2^*
r_2 =  z.$\smallskip

Proposition \ref{p direct sums} implies that $T ({A_{1}})$ and
$T({A_{2}}) $ are orthogonal weak$^*$ closed inner ideals of $B$,
whose direct sum is $B$. Thus, the unit of $B$ decomposes in the
form $1_{B} = v +w$, where $v\in T(A_1)$ and $w\in T(A_2)$. Since
$v$ and $w$ are orthogonal we have $1_{B} = 1_{B} 1_{B}^* =
(v+w)(v^*+w^*) = vv^* + ww^*$ and $1_B = v^*v + w^*w$, which shows
that $vv^*$ and $ww^*$ (respectively, $v^* v$ and $w^* w$) are two
orthogonal projections in $B$ whose sum is $1_B$. It follows that
$v v^* y = y v^* v = y$ for every element $y$ in
$T(A_1)$.\smallskip

It can be checked that $u = v+r_2$ is a unitary element in $B$ and
the mapping $$\Phi : (B,\circ_{u},\sharp_{u}) \to (B,\circ,*), \
x\mapsto x u^*$$ is a Jordan $^*$-isomorphism. By noticing that
$T(A_2)$ is a weak$^*$ closed inner ideal of $B$, $r_2 \in T(A_2)$
and $B_2(r_2)$ is the smallest weak$^*$ closed inner ideal
containing $r_2$, we have $T(A_2) = B_2 (r_2)$. Since $B$
decomposes in the form $$B_2 (v) \oplus^{\infty} B_2 (r_2) = B =
T(A_1) \oplus^{\infty} T(A_2),$$ we deduce that $v\in T(A_1)
\subseteq B_2 (v)$, which gives $T(A_1) = B_2 (v)$. We also have
$$B = \Phi (B_2 (v)) \oplus^{\infty} \Phi (B_2 (r_2)) = B_2 (vv^*)
\oplus^{\infty} B_2 (r_2 r_2^*)$$ $$ = vv^* B vv^* \oplus^{\infty}
r_2 r_2^* B r_2 r_2^*.$$

The mapping $\Phi|_{B_2 (r_2)} S_2$ is a Jordan $^*$-isomorphism
from $A_2$ onto $B_2 (r_2 r_2^*)$. Since $B_2 (r_2 r_2^*)$ is a
weak$^*$ closed ideal of $B$ and every Jordan $^*$-isomorphism
preserves the Murray-von Neumann decomposition we deduce that
$$\Phi|_{B_2 (r_2)} S_2 (A_2) \subseteq B_2 =
B_{II_1}\oplus^{\infty} B_{p\infty},$$ that is, $\Phi (B_2 (r_2))
\subseteq B_2$. We can similarly prove that $\Phi^{'} (B_2 (r_2))
\subseteq B_2$, where $\Phi^{'} (x) := u^* x$. Having in mind that
$B_2$ is a weak$^*$ closed ideal of $B$ we have $$T(A_2) = h_2
\circ_{r_2} S_2 (A_2) \subseteq \frac12 (h_1 r_2^* S(A_2) + S(A_2)
r_2^* h_2 )$$ $$= \frac12 (h_1 \Phi^{'} (B_2 (r_2)) + \Phi (B_2
(r_2)) h_2 )\subseteq \frac12 (h_1 B_2 + B_2 h_2 ) \subseteq
B_2.$$ It follows that $T^{-1} (B_2) \subseteq A_2$, and hence
$T(A_2) = B_2$. Thus
$$T(A_{1}) = T(A_{1}^{\perp\perp}) = T(A_{1}^{\perp})^{\perp} =
T(A_{2})^{\perp} = B_{2}^{\perp} = B_{1}.$$

Finally, since $A_1$ and $B_1$ are finite type I von Neumann
algebras, Proposition \ref{p finite type I} proves that
$T|_{A_{1}} : A_{1} \to B_{1}$ is continuous, which shows that $T$
enjoys the same property.
\end{proof}

\begin{remark}\label{r counterexample}{\rm  Let $A$ be a properly infinite 
or a (finite) type II$_1$ von Neumann algebra. Corollary \ref{prop-infinite} 
shows that every orthogonality preserving linear map from $A$ into a 
C$^*$-algebra is continuous.  We shall present an example showing that 
a similar statement doesn't hold when $A$ 
is replaced with a finite type I von Neumann algebra. In other words, 
the hypothesis of $T$ being surjective can not be removed in 
Theorem \ref{t vonNeumann}.\smallskip 

It is well know that a von Neumann algebra $A$ is type I and finite 
if and only if $A$ decomposes in the form $$A= \bigoplus_{i\in
I}^{\ell_{\infty}} C(\Omega_i, M_{m_i} (\mathbb{C})),$$ were the
$\Omega_i$'s are hyperstonean compact Hausdorff spaces and $(m_i)$
is a family of natural numbers  (c.f. \cite[Theorem V.1.27]{Tak}). 
In particular, every abelian von Neumann algebra is type I and finite.\smallskip

Let $K$ be an infinte (hyperstonean) compact set. By \cite[Example in page 142]{Jar}, there 
exists a discontinuous orthogonality preserving linear map 
$\varphi : C(K) \to \mathbb{C}$. Let $T: C(K) \to C(K) \oplus^{\infty} \mathbb{C}$ 
be the linear mapping defined by $T(f) := (f, \varphi(f))$ $(f\in C(K)).$ It is easy to 
check that $T$ is discontinuous and biorthogonality preserving but not surjective.   
}\end{remark}

Following \cite{ArauJar} and \cite{BurFerGarMarPe}, a linear map
$T$ between algebras $A$, $B$ is called \emph{separating} or
\emph{zero-product preserving} if $a b = 0$ implies $T (a) T(b)=
0$, for all $a, b$ in A; it is called \emph{biseparating} if
$T^{-1} :B\to A$ exists and is also separating. J. Araujo and K.
Jarosz conjectured in \cite[Conjecture 1]{ArauJar} that every
biseparating map between C$^*$-algebras is automatically
continuous. We can now give a partial positive answer to this
conjecture.\smallskip

Let $T: A \to B$ be a symmetric linear mapping between two
C$^*$-algebras. Suppose that $T$ is separating. Then for every
$a,b$ in $A$ with $a\perp b$, we have $T(a) T(b)^* = T(a) T(b^*)
=0$, because $T$ is separating. We can similarly prove that
$T(b)^* T(a) =0$, which shows that $T$ is orthogonality
preserving. The following result follows now as a consequence of
Theorem \ref{t vonNeumann}.

\begin{corollary}
\label{c biseparating symmetric} Let $T: A\to B$ be a biseparating
symmetric linear map between von Neumann algebras. Then $T$ is
continuous.$\hfill\Box$
\end{corollary}

\end{document}